\documentclass[conference]{IEEEtran} %!PN

\usepackage{latexsym}
\usepackage{amsfonts}
\usepackage{amsmath}
\usepackage{mathtools}
\usepackage{graphicx}
\usepackage{cite}
\usepackage{multicol}

\def\BibTeX{{\rm B\kern-.05em{\sc i\kern-.025em b}\kern-.08em
    T\kern-.1667em\lower.7ex\hbox{E}\kern-.125emX}}

\newtheorem{theorem}{Theorem}
\newtheorem{remark}{Remark}
\newtheorem{defin}{Definition}
\newtheorem{corollary}[theorem]{Corollary}
\newtheorem{lemma}[theorem]{Lemma}
\newtheorem{example}{Example}
\allowdisplaybreaks

\begin{document}

%\title{Using the Style File IEEEtran.sty} 
\title{Diffusive Representations\\ for the Numerical Evaluation of Fractional Integrals} %!PN

\author{\IEEEauthorblockN{Kai Diethelm}\IEEEauthorblockA{FANG, 
	Technical University of Applied Sciences W\"urzburg-Schweinfurt,
	Ignaz-Sch\"on-Str.\ 11, 97421 Schweinfurt, Germany\\ 
	Email: kai.diethelm@thws.de}}

\markboth{ICFDA'22 Proceedings}
{Diethelm: Diffusive Representations for Fractional Integrals} %!PN

\maketitle

\begin{abstract}
	Diffusive representations of fractional differential and integral operators
	can provide a convenient means to construct efficient numerical algorithms
	for their approximate evaluation. In the current literature, many different
	variants of such representations have been proposed. 
	Concentrating on Riemann-Liouville integrals whose order is in $(0,1)$,
	we here present a 
	general approach that comprises most of these variants as 
	special cases and that allows a detailed investigation of the 
	analytic properties of each variant. The availability of this information 
	allows to choose concrete numerical methods for handling the 
	representations that exploit the specific properties, thus allowing to
	construct very efficient overall methods.
\end{abstract}


\section{Introduction}
\label{sec:intro}

Owing to the non-local character of fractional order 
differential and integral operators, 
their numerical evaluation in their traditional
representation is a computationally much more complex task than
the analog evaluation of their integer order counterparts, 
both with respect to run time and with respect to memory requirements,
especially when an evaluation at many points is required.
Diffusive representations of the fractional operators \cite{Mo1998,Mo2005}, 
also known as infinite state representations \cite{HSL2019}, can be used as
a foundation upon which algorithms can be constructed 
that do not have such a high complexity; indeed in an asymptotic sense
(i.e.\ when the number $N$ of evaluation points is large)
these methods require only $O(N)$ operations and an $O(1)$ 
amount of memory, which is the same as observed in methods
for integer order problems. In this paper, we shall investigate
representations of this type for Riemann-Liouville integrals
\begin{equation}
	\label{eq:def-rli}
	J_a^\alpha f(t) 
	= \frac1 {\Gamma(\alpha)} 
			\int_a^t (t-\tau)^{\alpha-1} f(\tau) \, \mathrm d \tau
\end{equation}
of order $\alpha > 0$ with starting point $a \in \mathbb R$
for functions $f \in C[a,b]$. This continues the author's
recent investigations of similar properties of fractional differential operators 
\cite{Di2021,Di2022a,Di2023a}.
The main technical goal will be to provide the basic elements of a general
theory of such representations that comprises many previously suggested
concrete approaches as special cases. In this context, we concentrate on those 
properties of the diffusive representations that are particularly relevant when
one attempts to design efficient numerical methods.

\section{The general approach}

Our first result is a general fundamental representation for
Riemann-Liouville integrals of a given continuous function $f$. 
Here and in the following, given a number $\alpha \in (0, \infty) \setminus \mathbb N$
that denotes the order of the fractional operator under consideration,
we use the notational conventions
\begin{equation}
	\label{eq:def-n-c}
	n = \lceil \alpha \rceil
	\quad \mbox{ and } \quad
	c_\alpha = \frac{\sin \pi \alpha}\pi 
				\cdot \prod_{\ell=1}^{n-1} \frac 1 {\ell-\alpha}.
\end{equation}
In the case $\alpha \in (0,1)$ (which is the case that arises in the vast majority of
practical applications of fractional calculus based models), this reduces to
\begin{equation}
	\label{eq:def-n-c-1}
	n = 1
	\quad \mbox{ and } \quad
	c_\alpha = \frac{\sin \pi \alpha}{\pi}.
\end{equation}
Moreover, we will use the following definition.

\begin{defin}
	\label{def:admissible}
	A function $\psi : \Omega \to (0, \infty)$ is called an 
	\emph{admissible transformation} if it has the following properties:
	\begin{itemize}
	\item Its domain $\Omega$ is a non-empty open interval.
	\item $\psi \in C^1(\Omega)$.
	\item $\psi$ is strictly monotonically increasing.
	\item $\lim_{\omega \to \inf \Omega} \psi(\omega) = 0$ and 
		$\lim_{\omega \to \sup \Omega} \psi(\omega) = +\infty$.
	\end{itemize}
\end{defin}

With these conventions, our first main result reads as follows.

\begin{theorem}
	\label{thm:diff-rep-fund}
	Let $f \in C[a, b]$ with some real numbers $a < b$, and
	let $\alpha > 0$, $\alpha \notin \mathbb N$.
	Moreover, assume that $\psi : \Omega \to (0, \infty)$ 
	is an admissible transformation.
	Then, for every $t \in [a, b]$, the Riemann-Liouville integral
	of order $\alpha$ of the function $f$ can be expressed in the
	form of the \emph{diffusive representation}
	\begin{subequations}
		\label{eq:diff-rep-complete}
	\begin{equation}
		\label{eq:diff-rep}
		J_a^\alpha f(t) = \int_\Omega \phi(t, \omega) \, \mathrm d \omega
	\end{equation}
	with
	\begin{align}
		\label{eq:phi1}
		\phi(t, \omega) 
		&= c_\alpha \psi'(\omega) (\psi(\omega))^{n-\alpha-1} \\
		& \nonumber 
			\qquad \times
			\int_a^t (t - \tau)^{n-1} \exp(- (t - \tau) 
						\psi(\omega)) f(\tau) \, \mathrm d \tau.
	\end{align}
	\end{subequations}
\end{theorem}

The proof of this result will be given in Section \ref{sec:proof-fund}.

The essential basis for all the constructions of the numerical algorithms 
that exploit such diffusive representations is the observation that, 
for any fixed $\omega \in \Omega$,
the function $\phi(\cdot, \omega)$ can be characterized as being 
the unique solution to an initial value problem for a
very simple differential equation 
of \emph{integer} order. The precise structure of this differential equation is 
given in the following theorem. 
Like Theorem \ref{thm:diff-rep-fund}, we shall also prove this
result in Section \ref{sec:proof-fund}.

\begin{theorem}
	\label{thm:ivp}
	Assume the hypotheses of Theorem \ref{thm:diff-rep-fund}. 
	Then, for any $\omega \in \Omega$, 
	the function $\phi(\cdot, \omega)$ is the unique solution 
	on the interval $[a,b]$ to the $n$-th order
	differential equation
	\begin{subequations}
		\label{eq:ivp}
	\begin{eqnarray}
		\label{eq:ode}
		\lefteqn{\sum_{k=0}^n \binom{n}{k} (\psi(\omega))^{n-k} 
				\frac{\partial^k \phi}{\partial t^k} (t, \omega)} \\
		&=& c_\alpha \psi'(\omega) (\psi(\omega))^{n-1-\alpha} (n-1)! \cdot f(t) 
		\nonumber
	\end{eqnarray}
	subject to the initial conditions
	\begin{equation}
		\label{eq:ic}
		\frac{\partial^k \phi}{\partial t^k}(a, \omega) = 0
		\qquad (k = 0, 1, \ldots, n-1).
	\end{equation}
	\end{subequations}
\end{theorem}

\begin{remark}
	The differential equation \eqref{eq:ode} only contains 
	derivatives with respect to $t$,
	and all these derivatives are of integer order.
	Therefore, when only looking at this equation, 
	one may consider $t$ as the variable and
	all other quantities (in particular, $\omega$) 
	as (fixed) parameters. Then, it is natural 
	to consider \eqref{eq:ode} as an ordinary, not a partial, 
	differential equation. More precisely, 
	it is an inhomogeneous $n$-th order linear ordinary 
	differential equation with constant
	coefficients subject to the homogeneous initial 
	conditions \eqref{eq:ic}.
\end{remark}

In the special case $0 < \alpha < 1$, the statement of 
Theorem \ref{thm:ivp} reduces to the following result.

\begin{corollary}
	\label{cor:ivp-1}
	Assume the hypotheses of Theorem \ref{thm:ivp}. Moreover,
	let $\alpha \in (0,1)$. Then, for any $\omega \in \Omega$, 
	the function $\phi(\cdot, \omega)$ is the unique solution 
	on the interval $[a,b]$ to the first order
	differential equation
	\begin{subequations}
		\label{eq:ivp-1}
	\begin{equation}
		\label{eq:ode-1}
		\frac{\partial \phi}{\partial t} (t, \omega) 
		= - \psi(\omega) \phi(t, \omega) 
			+ c_\alpha \psi'(\omega) (\psi(\omega))^{-\alpha} f(t)
	\end{equation}
	subject to the initial condition
	\begin{equation}
		\label{eq:ic-1}
		\phi(a, \omega) = 0.
	\end{equation}
	\end{subequations}
\end{corollary}

\begin{remark}
	\label{rem:specialcases}
	A number of special cases of such representations 
	(i.e., special choices of the 
	transformation $\psi$ and the corresponding domain 
	$\Omega$) have been 
	discussed in the literature. We mention the following examples:
	\begin{itemize}
	\item $\Omega = (0, \infty)$ and $\psi(\omega) = \omega^2$ 
		\cite{YA2002},
	\item $\Omega = (0, \infty)$ and $\psi(\omega) = \omega^{1 - \alpha}$ 
		for $0 < \alpha < 1$ \cite{Ch2005,SC2006},
	\item $\Omega = (-\infty, \infty)$ and $\psi(\omega) = \mathrm e^\omega$ 
		\cite{Di2022a,Di2023a}.
	\end{itemize}
	Other options that do not seem to have been discussed so far
	might be, e.g.,
	\begin{itemize}
	\item $\Omega = (0, 1)$ and $\psi(\omega) = \tan (\omega \pi/2)$ and
	\item $\Omega = (0, 1)$ and $\psi(\omega) = \omega^{\sigma} / (1 - \omega)^{\rho}$
		with some $\sigma, \rho > 0$.
	\end{itemize}
\end{remark}

\section{Properties of the function $\phi$}
\label{sec:prop-phi}

The main motivation of this paper is to provide a basis upon which
one can construct
efficient numerical methods for the calculation of $J_a^\alpha f(t_\nu)$,
$\nu = 1, 2, \ldots, N$, with a potentially very large value of $N$.
From Theorems \ref{thm:diff-rep-fund} and \ref{thm:ivp}, we can conclude 
that a possible option in this context is to use the representation 
\eqref{eq:diff-rep} for $t = t_\nu$ and to evaluate this integral numerically with the
help of a suitably chosen quadrature formula. The integrand 
$\phi(t_\nu, \cdot)$ needs to be known of course, and this can
be accomplished by numerically solving the initial value problem \eqref{eq:ivp}
with an appropriate algorithm. When taking a concrete decision for the
specific choices of these two numerical schemes, it is important to take the
properties of the function $\phi$ into account. The numerical 
methods should be chosen in a way that is suitable for functions with 
these properties. Therefore, we now provide an analysis of the function $\phi$,
attempting to establish those properties that are relevant in the context
of choosing the numerical methods. It will turn out that the properties
under investigation strongly depend on which admissible transformation
$\psi$ has been selected.

\subsection{The integration of $\phi$}

We begin with some properties that are related to the integration 
in eq.~\eqref{eq:diff-rep}. Our first result in this context deals 
with the smoothness of the integrand in 
this equation, i.e.\ the function $\phi(t, \cdot)$ for a fixed value of $t$.
This is an important feature because, from classical results in the theory of 
numerical integration \cite{BP2011}, we know that integrals can be numerically
computed with high accuracy with a relatively small effort if the integrand is 
many times differentiable (with respect to the integration variable)
whereas this requires a much higher effort if the integrand has only a few derivatives.

\begin{theorem}
	Under the assumptions of Theorem \ref{thm:diff-rep-fund}, we have:
	\begin{itemize}
	\item If the transformation function $\psi$ satisfies $\psi \in C^k(\Omega)$
		with some $k \in \mathbb N$ then, for any $t \in [a, b]$,
		we have $\phi(t, \cdot) \in C^{k-1}(\Omega)$.
	\item If the transformation function $\psi$ satisfies 
		$\psi \in C^\infty(\Omega)$
		then, for any $t \in [a, b]$,
		we have $\phi(t, \cdot) \in C^\infty(\Omega)$.
	\end{itemize}
\end{theorem}

Bearing in mind the reasoning mentioned above, it is therefore advisable to 
use functions $\psi$ that are differentiable very frequently, ideally being
in $C^\infty(\Omega)$.

\emph{Proof.}
	This is an immediate consequence of the representation \eqref{eq:phi1}.
$\Box$

A second aspect that is relevant in the context of identifiying a numerical
integration method that performs well for the integrand in question is the
integrand's asymptotic behaviour as the integration variable tends towards
the ends of the integration interval. In this context, we can provide the 
following results. 
%Here, the notation $A(y) \sim B(y)$ for some quantities $A$
%and $B$ that depend on a common variable $y$ means that there exist
%two positive constants $C_1, C_2$ such that $C_1 \le |A(y) / B(y)| \le C_2$
%for all values of $y$ under consideration.

\begin{theorem}
	\label{thm:asymp}
	Assume the hypotheses of Theorem \ref{thm:diff-rep-fund}, and
	let $t \in [a, b]$ be fixed.
	\begin{enumerate}
	\item There exists some 
		constant $C \in \mathbb R$ such that
		\[
			|\phi(t, \omega)| \le C \psi'(\omega) (\psi(\omega))^{-\alpha-1}
			\quad \mbox{ for } \omega \to \sup \Omega.
		\]
%	\item If additionally $f(t) \ne 0$ and ??? then there exists some 
%		constant $C \in \mathbb R$ such that
%		\[
%			|\phi(t, \omega)| \ge C \psi'(\omega) (\psi(\omega))^{-\alpha}
%			\quad \mbox{ for } \omega \to \sup \Omega.
%		\]
	\item There exists some 
		constant $C \in \mathbb R$ such that
		\[
			|\phi(t, \omega)| \le C \psi'(\omega) (\psi(\omega))^{n-\alpha-1}
			\quad \mbox{ for } \omega \to \inf \Omega.
		\]
%	\item If additionally $f(t) \ne 0$ and ??? then there exists some 
%		constant $C \in \mathbb R$ such that
%		\[
%			|\phi(t, \omega)| \ge C \psi'(\omega) (\psi(\omega))^{-\alpha}
%			\quad \mbox{ for } \omega \to \inf \Omega.
%		\]
	\end{enumerate}
\end{theorem}

\begin{remark}
	For the interpretation of these bounds, it is useful to remember that,
	by definition of the admissible transformation $\psi$, we have
	\[
		\lim_{\omega \to \sup \Omega} \psi(\omega) = \infty
		\mbox{ and }
		\lim_{\omega \to \inf \Omega} \psi(\omega) = 0	
	\]
	(see the last item in the list of defining properties of the admissible
	transformations in Definition \ref{def:admissible}).
\end{remark}

\emph{Proof of Theorem \ref{thm:asymp}.}
	From \eqref{eq:phi1}, using the substitution 
	$\rho = (t-\tau) \psi(\omega)$, we find that
	\allowdisplaybreaks[0]
	\begin{align*}
		| \phi(t, \omega) |
		& \le |c_\alpha| \psi'(\omega) (\psi(\omega))^{-\alpha-1} 
				\sup_{u \in [a,b]} | f (u) | \\
		& \qquad \times
				\int_0^{(t-a) \psi(\omega)} \rho^{n-1} 
								\mathrm e^{-\rho} \, \mathrm d \rho.
	\end{align*}
	\allowdisplaybreaks%
	Depending on whether we are dealing with part 1 or part 2 of the
	statement, we now continue by estimating the integral on the right-hand
	side of this inequality in two different ways.
	
	For part 1, we note that
	\[
		\int_0^{(t-a) \psi(\omega)} \rho^{n-1} 
								\mathrm e^{-\rho} \, \mathrm d \rho
			< \int_0^\infty \rho^{n-1} 
								\mathrm e^{-\rho} \, \mathrm d \rho
			= (n-1)!
	\]
	which completes the proof in this case.
	
	To prove part 2, we use the fact that $\mathrm e^{-\rho} \le 1$ for all values
	of $\rho$ inside the integration range. This allows us to estimate
	\begin{align*}
		\int_0^{(t-a) \psi(\omega)} \rho^{n-1} 
								\mathrm e^{-\rho} \, \mathrm d \rho
		& \le \int_0^{(t-a) \psi(\omega)} \rho^{n-1}  \, \mathrm d \rho \\
		& = \frac 1 n \left( (t-a) \psi(\omega) \right)^n
		= C' (\psi(\omega))^n
	\end{align*}
	with some positive constant $C'$, and the proof is complete in this case too.
$\Box$

\begin{example}
	For $\psi(\omega) = \mathrm e^\omega$, we have $\Omega = (-\infty, \infty)$.
	Then, the function $\psi$ is in $C^\infty(\Omega)$ and 
	$\psi(t, \omega) = O( \mathrm e^{-\alpha \omega})$ for $\omega \to \infty$ 
	and $\psi(t, \omega) = O(\mathrm e^{(n-\alpha)\omega})$ for $\omega \to -\infty$. 
	Therefore, an integral with
	this integrand is highly suitable for being numerically computed by dividing the
	integration range up into $(-\infty, 0)$ and $(0, \infty)$, using a linear
	transformation of variables such that the new integrands on the subintervals 
	asymptotically behave as
	$c \mathrm e^{-\omega}$ and $c \mathrm e^{\omega}$, respectively, 
	and with approximating
	the integral over each part with the help of, e.g., a Gauss-Laguerre quadrature
	formula \cite{DR}. This special case has been discussed in detail 
	in \cite{Di2022a,Di2023a}.
\end{example}

\subsection{The solution of the initial value problem}

Next we address the aspects that are related to the
initial value problem \eqref{eq:ivp}. The main result here refers to the 
smoothness of the solution $\phi(\cdot, \omega)$ 
to this initial value problem as a function of
the variable according to which we differentiate in eq.~\eqref{eq:ode}:

\begin{theorem}
	\label{thm:ivp-smooth}
	Assume the hypotheses of Theorem \ref{thm:diff-rep-fund}.
	Moreover, let $f \in C^\ell[a, b]$ with some $\ell \in \mathbb N_0$.
	Then, for any $\omega \in \Omega$, 
	$\phi(\cdot, \omega) \in C^{\ell+n}[a,b]$.
\end{theorem}

We shall once again prove this result in Section \ref{sec:proof-fund}.

\begin{remark}
	\label{rem:dissip}
	Let us consider the case $0 < \alpha < 1$. Then the differential equation
	that $\phi$ solves takes the form \eqref{eq:ode-1}.
	Introducing the notation
	\[
		F_\omega(t, Y) 
		= - \psi(\omega) Y + c_\alpha \psi'(\omega) (\psi(\omega))^{-\alpha} f(t)
	\]
	so that the right-hand side of the differential equation 
	\eqref{eq:ode-1} is $F_\omega(t, \phi(t, \omega))$,
	we can see (taking into consideration that our assumptions 
	on the function $\psi$
	imply that $\psi(\omega) \ge 0$ for all $\omega$) that
	\[
		(F_\omega(t, Y_1) - F_\omega(t, Y_1)) (Y_1 - Y_2) 
		= - \psi(\omega) (Y_1 - Y_2)^2
		\le 0
	\]
	for all $Y_1, Y_2 \in \mathbb R$.
	Hence we conclude that the differential equation \eqref{eq:ode-1} is 
	dissipative \cite[Definition 8.58]{Plato}. This is a feature that is very
	helpful when trying to derive error bounds for (implicit) numerical
	solvers that do not massively overestimate the true errors \cite[Section 8.9]{Plato}.
\end{remark}

\begin{remark}
	In practical cases, the function $f$ to be fractionally integrated is the
	solution to a fractional differential equation or a function closely related 
	to that. It is well known \cite[Section 6.4]{Di2010} that such functions
	tend to be continuous but not differentiable at the initial point $a$,
	although differentiability at the other points of the interval $[a,b]$
	can usually be observed. This fact limits the practical applicability
	of Theorem \ref{thm:ivp-smooth}. However, the precise nature of the
	nondifferentiable components of the functions under consideration is
	often well understood. This information can in principle be exploited
	to obtain a more precise description of the smoothness properties 
	of the function $\phi$. We intend to address this question in future work.
\end{remark}

\section{Numerical Schemes}
\label{sec:num}

We now briefly indicate how the findings described above 
can be used to construct efficient numerical methods. To this end,
we shall in this section restrict our attention to the practically
most important case $0 < \alpha < 1$. 

The basic idea 
is to use, for a given function $f$ and a given data point $t_k$, 
a quadrature formula to evaluate the integral in 
eq.\ \eqref{eq:diff-rep}, i.e.\ we write
\begin{equation}
	\label{eq:qf}
	J_a^\alpha f(t_k) 
	\approx \sum_{m=1}^M w_{m,M} \tilde \phi(t_k, \omega_{m,M})
\end{equation}
with suitable weights $w_{m,M} \in \mathbb R$ and nodes 
$\omega_{m,M} \in \Omega$ ($m = 1, 2, \ldots, M$) that should be
chosen in a way that is adapted to the asymptotic behaviour of the
function $\phi$ indicated in Theorem \ref{thm:asymp}.
In eq.~\eqref{eq:qf}, the function $\tilde \phi$ ideally should be identical to
$\phi$ itself, but in practical applications this is likely to be unknown,
and therefore one usually needs to replace it by an approximation.

To obtain this approximation, one exploits the fact that $\phi$ solves the 
initial value problem \eqref{eq:ivp} and solves this problem numerically.
(Note that, in view of our assumption $0 < \alpha < 1$, this involves a first
order differential equation.)
For the sake of simplicity, we shall here concentrate on the case that this
is done with the help of a one-step solver like the backward Euler method.
This allows us to obtain the required solution at the point $t_k$ based
on the knowledge of the solution at the preceding point $t_{k-1}$
(which is known from the initial condition for $k=1$ and has been 
computed in the preceding time step otherwise). From eq.~\eqref{eq:qf},
we can see that we need to solve the initial value problem not just once,
but actually $M$ times, namely for the parameter values
$\omega = \omega_{m,M}$ ($m = 1, 2, \ldots, M$).

In these cases, the differential equation \eqref{eq:ode}
(which, in this case, has the form given in \eqref{eq:ode-1})
can be written as
\begin{equation}
	\label{eq:ode-1a}
	\frac{\partial \phi}{\partial t} (t, \omega_{m,M}) 
	= - \lambda_{m,M} \phi(t, \omega_{m,M}) + f_{m,M}(t)
\end{equation}
with 
$f_{m,M}(t) = c_\alpha \psi'(\omega_{m,M}) (\psi(\omega_{m,M}))^{-\alpha} f(t)$
and $\lambda_{m,M} = \psi(\omega_{m,M})$.
Taking into account that eq.~\eqref{eq:ode-1a} is dissipative
(see Remark \ref{rem:dissip}), the differential equation \eqref{eq:ode-1a} 
will be classified as stiff if 
\[
	\lambda_{m,M} = \psi(\omega_{m,M}) \gg 0,
\]
cf.~\cite[Section 8.9]{Plato}.

This observation explains the main challenge that one has to face when
constructing numerical methods based on this approach: 
If $\omega_{m,M}$ is such that $\psi(\omega_{m,M})$
is large then the differential equation \eqref{eq:ode-1a} becomes very stiff,
and this is likely to cause difficulties in the numerical solution.
In particular, reasonable results can only be expected when an implicit method
is used \cite[Section 8.9]{Plato}, e.g.\ the backward Euler method or the trapezoidal method.
Nevertheless, it might be advisable to avoid such points as much as possible. 
Clearly, in view of the fact that $\psi(\omega) \to \infty$ as $\omega \to
\sup \Omega$ by definition of an admissible transformation,
we cannot avoid such points completely in the limit case $M \to \infty$
since then we must not have any nonempty subintervals of $\Omega$
that do not contain any $\omega_{m,M}$. But in practical applications,
$M$ is finite, and then the discrete points 
$\{ \omega_{m,M} : m = 1, 2, \ldots, M \}$ leave gaps between them.
In this situation, since we want to have a good approximation to the
integral in eq.\ \eqref{eq:diff-rep}
with the formula \eqref{eq:qf}, one should aim to place those gaps 
(where the value of the integrand does not contribute to the
quadrature sum)
in regions where $|\phi(t, \omega)|$ is small (i.e.\ in regions which do not
contribute to the overall value of the integral in a very substantial way).
By Theorem \ref{thm:asymp}, this is the case
if $\psi'(\omega) (\psi(\omega))^{-\alpha-1}(\omega)$ is small.
%for $\omega$ near $\sup \Omega$ and
%$\psi'(\omega) (\psi(\omega))^{n-\alpha-1}(\omega)$ is small for
%$\omega$ near $\inf \Omega$. From the proposals indicated in 
%Remark \ref{rem:specialcases}, the special case
%\begin{equation}
%	\label{eq:psi-speziell}
%	\psi(\omega) = \frac{\omega^{1/(1-\alpha)}}{(1 - \omega)^{1/\alpha}}
%	\mbox{ with }
%	\Omega = (0, 1)
%\end{equation}
%that we have plotted in Fig.\ \ref{fig:psi} on the interval $\Omega$
%for three different values of $\alpha$
%has all the required properties in the case $0 < \alpha < 1$ that we discuss
%in this subsection. This can be seen in Fig.\ \ref{fig:psistrich} where we have
%plotted the function
%\begin{equation}
%	\label{eq:bound-omega}
%	\min \left \{ \psi'(\omega) \psi(\omega)^{-\alpha-1},
%				\psi'(\omega) \psi(\omega)^{n-\alpha-1} \right \},
%			\,
%			\omega \in \Omega = (0, 1),
%\end{equation}
%which, by Theorem \ref{thm:asymp}, is (up to a multiplicative constant)
%an upper bound for $|\phi(t, \omega)|$.
%(To handle the case $\alpha > 1$, a small modification
%is required, but this does not lead to major difficulties.)
%
%\begin{figure}[htb]
%	\centering
%	\includegraphics[width=0.9\columnwidth]{psi}
%	\caption{\label{fig:psi}Plots of $\psi(\omega)$ vs.\ $\omega$ 
%		according to \eqref{eq:psi-speziell} for three
%		different values of $\alpha$.}
%\end{figure}
%
%\begin{figure}[htb]
%	\centering
%	\includegraphics[width=0.9\columnwidth]{bound-omega}
%	\caption{\label{fig:psistrich}Plots of the bound \eqref{eq:bound-omega} 
%		vs.\ $\omega$ for three different values of $\alpha$.}
%\end{figure}

We intend to devote a future separate
paper to a thorough and systematic 
search for admissible transformations with this property,
for suitable quadrature formulas that exploit this behaviour,
and for (implicit) solvers for initial value problems that can
be well combined with these quadrature formulas.

\section{Proofs of the main results}
\label{sec:proof-fund}

In this section, we will provide the proofs of our main
results that have not been given already, viz.\ Theorems 
\ref{thm:diff-rep-fund}, \ref{thm:ivp} and \ref{thm:ivp-smooth}.

\subsection{Proof of Theorem \protect{\ref{thm:diff-rep-fund}}}

For the proof of Theorem \ref{thm:diff-rep-fund}, we need 
some well known fundamental 
properties of the Gamma function that we recall here, 
namely its definition \cite[eq.~(2.1)]{Artin}
\begin{equation}
	\label{eq:def-gamma}
	\Gamma(z) = \int_0^\infty u^{z-1} \mathrm e^{-u} \, \mathrm d u
	\quad \mbox{for } z > 0,
\end{equation}
the functional equation \cite[eq.~(2.2)]{Artin}
\begin{equation}
	\label{eq:func-eq-gamma}
	z \Gamma(z) = \Gamma(z+1)
	\quad \mbox{for } z > 0,
\end{equation}
and the reflection formula \cite[eq.~(4.5)]{Artin}
\begin{equation}
	\label{eq:refl-gamma}
	\frac 1 {\Gamma(z)} = \frac{\sin \pi z}{\pi} \Gamma(1-z)
	\quad \mbox{for } z \in \mathbb R \setminus \mathbb N.
\end{equation}

In the first step of the proof,
we combine the definition \eqref{eq:def-rli} of $J_a^\alpha f(t)$ with
the reflection formula \eqref{eq:refl-gamma} and obtain
\[
	J_a^\alpha f(t) 
	= \frac{\sin \pi \alpha}\pi \Gamma(1-\alpha)
		 \int_a^t (t-\tau)^{\alpha-1} f(\tau) \, \mathrm d \tau.
\]
In the next step, we apply the functional equation \eqref{eq:func-eq-gamma}
$n-1$ times which, in view of the definition of $c_\alpha$ 
given in \eqref{eq:def-n-c}, yields
\[
	J_a^\alpha f(t) 
	= c_\alpha
		 \int_a^t \Gamma(n-\alpha) (t-\tau)^{\alpha-1} f(\tau) \, \mathrm d \tau.
\]
Noting that, by the definition \eqref{eq:def-n-c} of $n$, 
the argument $n - \alpha$ of the Gamma function in this equation 
is in $(0,1)$, we replace the Gamma function by the integral on the
right-hand side of \eqref{eq:def-gamma}, thus arriving at
\begin{align*}
	J_a^\alpha f(t) 
	&= c_\alpha
		 \int_a^t \int_0^\infty u^{n-\alpha-1} \mathrm e^{-u} \, \mathrm d u 
		 		\, (t-\tau)^{\alpha-1} f(\tau) \, \mathrm d \tau \\
	&= c_\alpha
		 \int_a^t \int_0^\infty \left( \frac{u}{t-\tau}\right)^{1-\alpha} \mathrm e^{-u} 
		 		u^{n-2} f(\tau) \, \mathrm d u \, \mathrm d \tau .
\end{align*}
The substitution $u = (t - \tau) \psi(\omega)$ in the inner integral then leads us
to
\begin{align*}
	J_a^\alpha f(t) 
	&= c_\alpha
		 \int_a^t \int_\Omega 
		 	(\psi(\omega))^{n-1-\alpha} \mathrm e^{-(t-\tau) \psi(\omega)} \\
	& \qquad \qquad \qquad \quad 
		\times (t - \tau)^{n-1}
		 	f(\tau) \psi'(\omega) \, \mathrm d \omega \, \mathrm d \tau.
\end{align*}
Under the given conditions, Fubini's Theorem allows to interchange the order of
integration, and this gives the desired representation 
\eqref{eq:diff-rep-complete}.~$\Box$

\subsection{An auxiliary result}

Next, we state an auxiliary result that we will need later in the proof
of Theorem \ref{thm:ivp}.

\begin{lemma}
	\label{lem:bin-coeff}
	Let $n \in \mathbb N$ and $\mu \in \{0, 1, 2, \ldots, n-1\}$.
	Then,
	\[
		\sum_{k=\mu}^{n} \binom{n}{k}  \binom{k}{\mu} (-1)^{k-\mu} = 0.
	\]
\end{lemma}

\emph{Proof.}
	\begin{align*}
		\MoveEqLeft{\sum_{k=\mu}^{n} \binom{n}{k} \binom{k}{\mu} (-1)^{k-\mu} } \\
		& = \sum_{k=\mu}^{n} \frac{n!}{k! (n-k)!} 
						\cdot \frac{k!}{\mu! (k-\mu)!} (-1)^{k-\mu}  \\
		& = \frac{n!}{\mu! (n-\mu)!} 
				\sum_{k=\mu}^{n} \frac{(-1)^{k-\mu} (n-\mu)!}{(k-\mu)! (n-k)!} \\
		& = \binom{n}{\mu} 
				\sum_{k=\mu}^{n} (-1)^{k-\mu} \binom{n - \mu}{n - k} \\
		& =   (-1)^{n -\mu} \binom{n}{\mu} 
				\sum_{\ell=0}^{n-\mu} (-1)^{\ell} \binom{n - \mu}{\ell} = 0
	\end{align*}
where, in the last step, we have used a well known identity for 
binomial coefficients \cite[eq.~3.1.7]{AS}.~$\Box$

\subsection{Proof of Theorem \protect{\ref{thm:ivp}}}

We start by collecting some basic properties of the function~$\phi$.
To this end, we shall use the notation
\[
	g_1(z) = z^{n-1}, \quad 
	g_2(z) = \mathrm e^{-z \psi(\omega)}
	\mbox{ and }
	g(z) = g_1(z) g_2(z).
\]
A $k$-fold differentiation of $g$ where $k \in \{ 0, 1, 2, \ldots, n \}$ leads
to
\begin{equation}
	\label{eq:g-k-diff}
	g^{(k)}(z) 
	= \sum_{\mu=0}^k \binom{k}{\mu} g_1^{(\mu)}(z) g_2^{(k-\mu)}(z).
\end{equation}
Since $g_1$ has an ($n-1$)-fold zero at $0$, it is clear that $g_1^{(\mu)}(0) = 0$
for $\mu = 0, 1, \ldots, n-2$.
We thus conclude 
\begin{equation}
	\label{eq:g0}
	g^{(k)}(0) = 0
	\mbox{ for } k = 0, 1, \ldots, n-2.
\end{equation}

By definition of $g$, it follows from \eqref{eq:phi1} that
\begin{equation}
	\label{eq:phi-g}
	\phi(t, \omega) 
	= c_\alpha^* \int_a^t g(t - \tau) f(\tau) \, \mathrm d \tau
\end{equation}
where 
\[
	c_\alpha^* = c_\alpha \psi'(\omega) (\psi(\omega))^{n-\alpha-1}.
\]
Hence, using \eqref{eq:g0}, a $k$-fold differentiation of \eqref{eq:phi-g} 
yields
\begin{equation}
	\label{eq:phi-k-diff}
	\frac{\partial^k \phi}{\partial t^k}(t, \omega) 
		= c_\alpha^*  \int_a^t g^{(k)}(t-\tau) f(\tau) \, \mathrm d \tau
	\mbox{  for } k = 0, 1, \ldots, n-1.
\end{equation}
Evaluating these expressions and eq.~\eqref{eq:phi-g} at $t=a$, 
we can see that $\phi$ satisfies the initial condition \eqref{eq:ic}.

Using the case $k=n-1$ of eq.~\eqref{eq:phi-k-diff} and differentiating once
more, we find
\begin{equation}
	\label{eq:phi-n-diff}
	\frac{\partial^n \phi}{\partial t^n}(t, \omega) 
		= c_\alpha^* \left(   g^{(n-1)}(0) f(t) 
							+ \int_a^t g^{(n)}(t-\tau) f(\tau) \, \mathrm d \tau	
					 \right).
\end{equation}
Taking into account eqs.\ \eqref{eq:g-k-diff} and \eqref{eq:g0}, we can see that
\[
	g^{(n-1)}(0) = g_1^{(n-1)}(0) g_2(0) 
		= \left[ \frac{\mathrm d^{n-1}}{\mathrm d z^{n-1}} z^{n-1} \right]_{z=0}
		= (n-1)!.
\]
Plugging this relation into \eqref{eq:phi-n-diff}, we obtain
\begin{equation}
	\label{eq:phi-n-diff2}
	\frac{\partial^n \phi}{\partial t^n}(t, \omega) 
		= c_\alpha^* \left( (n-1)! \cdot f(t) 
							+ \int_a^t g^{(n)}(t-\tau) f(\tau) \, \mathrm d \tau	
					\right).
\end{equation}

Exploiting eqs.~\eqref{eq:phi-k-diff} and \eqref{eq:phi-n-diff2}, we thus derive
\begin{eqnarray*}
	\lefteqn{\sum_{k=0}^n \binom{n}{k} (\psi(\omega))^{n-k} 
						\frac{\partial^k \phi}{\partial t^k}(t, \omega) } \\
	& = &  c_\alpha^* \Bigg( (n-1)! \cdot f(t)  \\
	&& \qquad
					+ \sum_{k=0}^n \binom{n}{k} (\psi(\omega))^{n-k} 
							\int_a^t g^{(k)}(t-\tau) f(\tau) \, \mathrm d \tau	
					\Bigg) \\
	& = &  c_\alpha^* \Bigg( (n-1)! \cdot f(t)  \\
	&& \qquad
		+ \int_a^t f(\tau) \sum_{k=0}^n \binom{n}{k} 
			(\psi(\omega))^{n-k} g^{(k)}(t-\tau) \, \mathrm d \tau	
					\Bigg) 
\end{eqnarray*}
In order to complete the proof of our desired result \eqref{eq:ode},
it now suffices to show that the sum inside of the integral on the right-hand
side of this equation vanishes identically. To see this, we write
\begin{eqnarray*}
	\lefteqn{\sum_{k=0}^n \binom{n}{k} (\psi(\omega))^{n-k} g^{(k)}(z)} \\
	&=& \sum_{k=0}^n \binom{n}{k} (\psi(\omega))^{n-k} 
		\sum_{\mu=0}^k \binom{k}{\mu} g_1^{(\mu)}(z) g_2^{(k-\mu)}(z) .
\end{eqnarray*}
In view of the definitions of $g_1$ and $g_2$, the derivatives of these functions 
are easy to compute, namely
\begin{align*}
	g_1^{(\mu)}(z) &= \begin{cases}
					\frac{(n-1)!}{(n-1-\mu)!} z^{n-1-\mu} 
						 & \mbox{ for } \mu = 0, 1, \ldots, n-1, \\
					0 & \mbox{ for } \mu = n,
				\end{cases} \\ 
	g_2^{(k-\mu)}(z) &= (-\psi(\omega))^{k-\mu} \mathrm e^{-z \psi(\omega)}.
\end{align*}
Thus,
\begin{eqnarray*}
	\lefteqn{\sum_{k=0}^n \binom{n}{k} (\psi(\omega))^{n-k} g^{(k)}(z)} \\
	&=& \sum_{k=0}^n \binom{n}{k} (\psi(\omega))^{n-k} 
							 \mathrm e^{-z \psi(\omega)} \\
	&& \quad \times
		\sum_{\mu=0}^{\min(k, n-1)} \binom{k}{\mu} 
			\frac{(n-1)!}{(n-1-\mu)!} z^{n-1-\mu} 
			(-\psi(\omega))^{k-\mu} \\
	&=&  \mathrm e^{-z \psi(\omega)} (n-1)! \sum_{\mu=0}^{n-1} 
		\frac{z^{n-1-\mu} (\psi(\omega))^{n-\mu}}{(n-1-\mu)!} \\
	&& \qquad\qquad\qquad\qquad\qquad \times
			\sum_{k=\mu}^{n}	\binom{n}{k}  \binom{k}{\mu} (-1)^{k-\mu}.
\end{eqnarray*}
The innermost sum is zero because of Lemma \ref{lem:bin-coeff}, and hence the 
entire expression is zero as required.~$\Box$

\subsection{Proof of Theorem \protect{\ref{thm:ivp-smooth}}}

For the sake of exposition, we provide two
 different proofs of Theorem \ref{thm:ivp-smooth}.

Firstly, using the notation $\Phi(t) = \phi(t, \omega)$ for some
arbitrary but fixed $\omega \in \Omega$, the differential equation 
\eqref{eq:ode} can be written in the form 
\begin{equation}
	\label{eq:ode-G}
	\Phi^{(n)}(t) = G(t, \Phi(t), \Phi'(t), \ldots, \Phi^{(n-1)}(t))		
\end{equation}
with
\begin{align*}
	\MoveEqLeft {G(t, z_0, z_1, \ldots, z_{n-1}) } \\
	& = - \sum_{k=0}^{n-1} \binom{n}{k} (\psi(\omega))^{n-k} z_k \\
		& \qquad \qquad 
			+ c_\alpha \psi'(\omega) (\psi(\omega))^{n-1-\alpha} (n-1)! \cdot f(t).
\end{align*}
Clearly, under the given assumptions, the function $G$ on the
right-hand side of eq.~\eqref{eq:ode-G} possesses $\ell$ continuous
partial derivatives with respect to $t$. The claim is then an immediate
consequence of some well known general properties of differential 
equations \cite[Chapter 1, Theorem 1.2]{CL}.~$\Box$

Alternatively, we can also resort to the representation of the $n$-th
partial derivative of $\phi(t, \omega)$ with respect to $t$ that was given 
in eq.~\eqref{eq:phi-n-diff}. It is evident that, under our assumptions, the
expression on the right-hand side of \eqref{eq:phi-n-diff} possesses $\ell$
continuous derivatives with respect to $t$, and this also implies our claim.
$\Box$

\section{Conclusions}

We have introduced a general framework for designing many possible
special cases of diffusive representations for fractional integral operators. 
Our abstract analysis of these
repesentations allowed us to establish a number of their properties. The
knowledge of these properties can be exploited to design fast and accurate
algorithms for the numerical evaluation of Riemann-Liouville integrals at
multiple data points. 

In view of the well known relationships between Riemann-Liouville integrals
and Caputo derivatives, the results can also be transferred to numerical
methods for the latter. Moreover, algorithms based on our findings can form 
a basic building block for fast solvers for fractional differential equations.

\addtolength{\labelsep}{-1mm}
\addtolength{\columnwidth}{-1mm}
\bibliographystyle{IEEEtran}
\bibliography{icsc}

\begin{thebibliography}{10}
\providecommand{\url}[1]{#1}
\csname url@samestyle\endcsname
\providecommand{\newblock}{\relax}
\providecommand{\bibinfo}[2]{#2}
\providecommand{\BIBentrySTDinterwordspacing}{\spaceskip=0pt\relax}
\providecommand{\BIBentryALTinterwordstretchfactor}{4}
\providecommand{\BIBentryALTinterwordspacing}{\spaceskip=\fontdimen2\font plus
\BIBentryALTinterwordstretchfactor\fontdimen3\font minus
  \fontdimen4\font\relax}
\providecommand{\BIBforeignlanguage}[2]{{%
\expandafter\ifx\csname l@#1\endcsname\relax
\typeout{** WARNING: IEEEtran.bst: No hyphenation pattern has been}%
\typeout{** loaded for the language `#1'. Using the pattern for}%
\typeout{** the default language instead.}%
\else
\language=\csname l@#1\endcsname
\fi
#2}}
\providecommand{\BIBdecl}{\relax}
\BIBdecl

\bibitem{Mo1998}
G.~Montseny, ``Diffusive representation of pseudo-differential
  time-operators,'' \emph{ESAIM, Proc.}, vol.~5, pp. 159--175, 1998.

\bibitem{Mo2005}
------, \emph{Repr\'esentation diffusive}.\hskip 1em plus 0.5em minus
  0.4em\relax Paris: Hermes Science/Lavoisier, 2005.

\bibitem{HSL2019}
M.~Hinze, A.~Schmidt, and R.~I. Leine, ``Numerical solution of fractional-order
  ordinary differential equations using the reformulated infinite state
  representation,'' \emph{Fract. Calc. Appl. Anal.}, vol.~22, pp. 1321--1350,
  2019.

\bibitem{Di2021}
K.~Diethelm, ``Fast solution methods for fractional differential equations in
  the modeling of viscoelastic materials,'' in \emph{Proc. 9th International
  Conference on Systems and Control}.\hskip 1em plus 0.5em minus 0.4em\relax
  Piscataway: IEEE, 2021, pp. 455--460.

\bibitem{Di2022a}
------, ``A new diffusive representation for fractional derivatives, part {II}:
  Convergence analysis of the numerical scheme,'' \emph{Mathematics}, vol.~10,
  p. 1245, 2022.

\bibitem{Di2023a}
------, ``A new diffusive representation for fractional derivatives, part {I}:
  Construction, implementation and numerical examples,'' in \emph{Fractional
  Differential Equations: Modeling, Discretization, and Numerical Solvers},
  A.~Cardone, M.~Donatelli, F.~Durastante, R.~Garrappa, M.~Mazza, and
  M.~Popolizio, Eds.\hskip 1em plus 0.5em minus 0.4em\relax Singapore: Springer
  Nature, 2023.

\bibitem{YA2002}
L.~Yuan and O.~P. Agrawal, ``A numerical scheme for dynamic systems containing
  fractional derivatives,'' \emph{J. Vibration Acoustics}, vol. 124, pp.
  321--324, 2002.

\bibitem{Ch2005}
A.~Chatterjee, ``Statistical origins of fractional derivatives in
  viscoelasticity,'' \emph{J. Sound Vibration}, vol. 284, pp. 1239--1245, 2005.

\bibitem{SC2006}
S.~J. Singh and A.~Chatterjee, ``Galerkin projections and finite elements for
  fractional order derivatives,'' \emph{Nonlinear Dynamics}, vol.~45, pp.
  183--206, 2006.

\bibitem{BP2011}
H.~Bra\ss\ and K.~Petras, \emph{Quadrature Theory}.\hskip 1em plus 0.5em minus
  0.4em\relax Providence: Amer. Math. Soc., 2011.

\bibitem{DR}
P.~J. Davis and P.~Rabinowitz, \emph{Methods of Numerical Integration},
  2nd~ed.\hskip 1em plus 0.5em minus 0.4em\relax San Diego: Academic Press,
  1984.

\bibitem{Plato}
R.~Plato, \emph{Concise Numerical Mathematics}.\hskip 1em plus 0.5em minus
  0.4em\relax Providence: Amer. Math. Soc., 2003.

\bibitem{Di2010}
K.~Diethelm, \emph{The Analysis of Fractional Differential Equations}.\hskip
  1em plus 0.5em minus 0.4em\relax Berlin: Springer, 2010.

\bibitem{Artin}
E.~Artin, \emph{The Gamma Function}.\hskip 1em plus 0.5em minus 0.4em\relax New
  York: Holt, Rinehart \& Winston, 1964.

\bibitem{AS}
M.~Abramowitz and I.~A. Stegun, \emph{Handbook of Mathematical
  Functions}.\hskip 1em plus 0.5em minus 0.4em\relax Washington: National
  Bureau of Standards, 1972.

\bibitem{CL}
E.~A. Coddington and N.~Levinson, \emph{Theory of Ordinary Differential
  Equations}.\hskip 1em plus 0.5em minus 0.4em\relax New York: McGraw-Hill,
  1955.

\end{thebibliography}

\end{document}